\documentclass[11pt]{amsart}
\usepackage{amsfonts}
\usepackage{amsmath}
\usepackage{amssymb}
\usepackage{color}
\usepackage[all]{xy}           
\usepackage{graphicx}
\usepackage{xspace}
\usepackage{axodraw}
\usepackage{tikz}
\def\shuff#1#2{\mathbin{
      \hbox{\vbox{\hbox{\vrule \hskip#2 \vrule height#1 width 0pt}\hrule}\vbox{\hbox{\vrule \hskip#2 \vrule height#1 width 0pt\vrule }\hrule}}}}
\def\shuffl{{\mathchoice{\shuff{5pt}{3.5pt}}{\shuff{5pt}{3.5pt}}{\shuff{3pt}{2.6pt}}{\shuff{3pt}{2.6pt}}}}

\def\shuffle{{\, \shuffl \,}}

\newtheorem{thm}{Theorem}

\newtheorem{lem}[thm]{Lemma}
\newtheorem{prop}[thm]{Proposition}
\newtheorem{defn}{Definition}

\newcommand*{\id}{{\mathrm{id}}}

\newcommand*{\Kb}{{\mathbb K}}

\newcommand*{\Nb}{{\mathbb N}}

\newcommand*{\un}{{\mathbf 1}}

 \def\arbrea{
 \begin{tikzpicture}
\draw[thick] (0,0) -- (0,-0.5);
                	\filldraw[green] (0,0) circle(2.5pt);
                	\filldraw[black] (0,-0.5) circle(2pt);
\end{tikzpicture}
}

 \def\arbreb{
 \begin{tikzpicture}
\draw[thick] (0,0) -- (0.5,-0.5);
\draw[thick] (0,0) -- (-0.5,-0.5);
                	\filldraw[black] (0.5,-0.5) circle(2pt);
                	\filldraw[black] (-0.5,-0.5) circle(2pt);
                	\filldraw[green] (0,0) circle(2.5pt);
\end{tikzpicture}
}

 \def\arbrec{
 \begin{tikzpicture}
\draw[thick] (0,0) -- (0,-0.5);
\draw[thick] (0,-0.5) -- (0,-1);
                	\filldraw[black] (0,-0.5) circle(2pt);
                	\filldraw[black] (0,-1) circle(2pt);
                	\filldraw[green] (0,0) circle(2.5pt);
\end{tikzpicture}
}

 \def\arbree{
 \begin{tikzpicture}
\draw[thick] (0,0) -- (0.5,-0.5);
\draw[thick] (0,0) -- (-0.5,-0.5);
\draw[thick] (0,0) -- (0,-0.5);
                	\filldraw[black] (0,-0.5) circle(2pt);
                	\filldraw[black] (0.5,-0.5) circle(2pt);
                	\filldraw[black] (-0.5,-0.5) circle(2pt);
                	\filldraw[green] (0,0) circle(2.5pt);
\end{tikzpicture}
}

 \def\arbref{
 \begin{tikzpicture}
\draw[thick] (0,0) -- (0.5,-0.5);
\draw[thick] (0,0) -- (-0.5,-0.5);
\draw[thick] (-0.5,-0.5) -- (-0.5,-1);
                	\filldraw[black] (0.5,-0.5) circle(2pt);
                	\filldraw[black] (-0.5,-0.5) circle(2pt);
                	\filldraw[black] (-0.5,-1) circle(2pt);
                	\filldraw[green] (0,0) circle(2.5pt);
\end{tikzpicture}
}

 \def\arbreg{
 \begin{tikzpicture}
\draw[thick] (0,0) -- (0,-0.5);
\draw[thick] (0,-0.5) -- (0.5,-1);
\draw[thick] (0,-0.5) -- (-0.5,-1);
                	\filldraw[green] (0,0) circle(2.5pt);
                	\filldraw[black] (0,-0.5) circle(2pt);
                	\filldraw[black] (0.5,-1) circle(2pt);
                	\filldraw[black] (-0.5,-1) circle(2pt);
\end{tikzpicture}
}

\def\baarr{
 \begin{tikzpicture}
\draw[ultra thick] (0.1,-0.4) -- (0.1,0); 
 \end{tikzpicture}}
 
\def\nci{
 \begin{tikzpicture}
\draw[thick] (0,-1) -- (0,0); 
 \end{tikzpicture}}

\def\ncii{
 \begin{tikzpicture}
\draw[thick]  (0,0) -- (1,0);
\draw[thick] (0,-1) -- (0,0); 
\draw[thick]  (1,-1) -- (1,0); 
 \end{tikzpicture}}

\def\nciii{
 \begin{tikzpicture}
\draw[thick]  (0,0) -- (2,0);
\draw[thick] (0,-1) -- (0,0); 
\draw[thick]  (1,-1) -- (1,0); 
\draw[thick]  (2,-1) -- (2,0); 
 \end{tikzpicture}}

\def\nciiib{
 \begin{tikzpicture}
\draw[thick]  (0,0) -- (1,0);
\draw[thick] (0,-1) -- (0,0); 
\draw[thick]  (0.5,-1) -- (0.5,-0.2); 
\draw[thick] (1,-1) -- (1,0); 
\draw[thick]  (0.5,-0.2) -- (1.5,-0.2);
\draw[thick]  (1.5,-1) -- (1.5,-0.2); 
 \end{tikzpicture}}

\def\nciiic{
 \begin{tikzpicture}
\draw[thick]  (0,0) -- (2,0);
\draw[thick] (0,-1) -- (0,0); 
\draw[thick]  (0.5,-1) -- (0.5,-0.2); 
\draw[thick]  (0.5,-0.2) -- (1.5,-0.2);
\draw[thick]  (1.5,-1) -- (1.5,-0.2); 
\draw[thick] (2,-1) -- (2,0); 
 \end{tikzpicture}}

\def\nciiia{
 \begin{tikzpicture}
\draw[thick]  (0,0) -- (2,0);
\draw[thick] (0,-1) -- (0,0); 
\draw[thick]  (1,-1) -- (1,-0.2); 
\draw[thick]  (2,-1) -- (2,0); 
 \end{tikzpicture}}

\def\nciiiia{
 \begin{tikzpicture}
\draw[thick]  (0,0) -- (3,0);
\draw[thick] (0,-1) -- (0,0); 
\draw[thick]  (0.75,-1) -- (0.75,-0.2); 
\draw[thick]  (1.5,-1) -- (1.5,-0.2); 
\draw[thick]  (2.25,-1) -- (2.25,-0.2); 
\draw[thick]  (3,-1) -- (3,0); 
 \end{tikzpicture}}

\def\nciiiiia{
 \begin{tikzpicture}
\draw[thick]  (0,0) -- (3,0);
\draw[thick] (0,-1) -- (0,0); 
\draw[thick]  (0.75,-1) -- (0.75,-.2); 
\draw[thick]  (1.5,-.2) -- (2.5,-.2); 
\draw[thick]  (1.5,-1) -- (1.5,-.2); 
\draw[thick]  (2.5,-1) -- (2.5,-.2); 
\draw[thick]  (3,-1) -- (3,0); 
 \end{tikzpicture}}

\def\nciiiiii{
 \begin{tikzpicture}
\draw[thick]  (0,0) -- (3,0);
\draw[thick] (0,-1) -- (0,0); 
\draw[thick]  (0.75,-1) -- (0.75,-.2); 
\draw[thick]  (1.5,-.2) -- (2.5,-.2); 
\draw[thick]  (1.5,-1) -- (1.5,-.2); 
\draw[thick]  (2,-1) -- (2,-.4); 
\draw[thick]  (2.5,-1) -- (2.5,-.2); 
\draw[thick]  (3,-1) -- (3,0); 
 \end{tikzpicture}}

\def\nciiiiiia{
 \begin{tikzpicture}
\draw[thick]  (0,0) -- (4,0);
\draw[thick] (0,-1) -- (0,0); 
\draw[thick]  (0.5,-1) -- (0.5,-.2); 
\draw[thick]  (1,-1) -- (1,0); 
\draw[thick]  (1.5,-.2) -- (3.5,-.2); 
\draw[thick]  (1.5,-1) -- (1.5,-.2); 
\draw[thick]  (3.5,-1) -- (3.5,-.2); 
\draw[thick]  (2,-1) -- (2,-.4); 
\draw[thick]  (2.5,-1) -- (2.5,-.2); 
\draw[thick]  (4,-1) -- (4,0); 
\draw[thick]  (4.5,0) -- (5.5,0); 
\draw[thick]  (4.5,-1) -- (4.5,0); 
\draw[thick]  (5.5,-1) -- (5.5,0); 

 \end{tikzpicture}}

\begin{document}

\title{The splitting process in free probability theory}

\author{Kurusch Ebrahimi-Fard}
\address{ICMAT,
		C/Nicol\'as Cabrera, no.~13-15, 28049 Madrid, Spain.
		On leave from UHA, Mulhouse, France}
         \email{kurusch@icmat.es}         
         \urladdr{www.icmat.es/kurusch}

\author{Fr\'ed\'eric Patras}
\address{Laboratoire J.-A.~Dieudonn\'e
         		UMR 6621, CNRS,
         		Parc Valrose,
         		06108 Nice Cedex 02, France.}
\email{patras@math.unice.fr}
\urladdr{www-math.unice.fr/$\sim$patras}

\voffset=10ex

\begin{abstract}
Free cumulants were introduced by Speicher as a proper analog of classical cumulants in Voiculescu's theory of free probability. The relation between free moments and free cumulants is usually described in terms of M\"obius calculus over the lattice of non-crossing partitions. In this work we explore another approach to free cumulants and to their combinatorial study using a combinatorial Hopf algebra structure on the linear span of non-crossing partitions. The generating series of free moments is seen as a character on this Hopf algebra. It is characterized by solving a linear fixed point equation that relates it to the generating series of free cumulants. These phenomena are explained through a process similar to (though different from) the arborification process familiar in the theory of dynamical systems, and originating in Cayley's work.
 \end{abstract}


\date{February 9, 2015}

\maketitle

\keywords{Keywords: Cumulants, free cumulants, non-crossing partitions, half-shuffles, half-unshuffles.}



\section{Introduction}
\label{sect:intro}

Free cumulants were introduced by Speicher as the proper analog of classical cumulants in Voiculescu's theory of free probability. The relation between free moments and free cumulants is usually described in terms of M\"obius calculus over the lattice of non-crossing partitions. See \cite{NiSp,Speicher,speiLoth,DVoi} for details. In this work we explore another approach to free cumulants and moments based on a  combinatorial Hopf algebra structure defined directly on non-crossing partitions. Starting with the nice work of Mastnak and Nica \cite{MastnakNica} on the logarithm of the $S$-transform, several approaches have been made in the past to combine the notion of free probability with that of -- combinatorial -- Hopf algebras \cite{BeBoLeSp,FriedMcKay1,FriedMcKay2}. The key ingredient, that distinguishes our approach from the classical understanding of the relation between free cumulants and moments, is the use of half-unshuffles. The latter are dual to half-shuffles, which were first introduced in the mid-1950's by Eilenberg and MacLane \cite{EM} with the aim to provide an axiomatic and algebraic understanding of products in the cohomology of topological spaces such as the $K(\pi,n)$, i.e., spaces with only one non trivial homotopy group, in dimension $n$.

In a previous article \cite{EPnc}, we showed that half-unshuffles defined on the double tensor algebra over a  non-commutative probability space $(A,\phi)$, where $A$ is a complex algebra with unit $1_A$ and $\phi$ is a $\mathbb{C}$-valued linear form on ${A}$, such that $\phi(1_{A}) = 1$, provide a rather simple way, to define and understand the connections between free moments and free cumulants. In fact, half-shuffles yield also a new approach to the classical notion of cumulants, that can be understood from this point of view as a plain ``commutative version'' of free cumulants. 

In the present article, we focus on subtler combinatorial properties of the calculus of free probabilities involving non-crossing partitions. 
Recall that non-crossing partitions are in bijection with planar rooted trees, so that any construction involving one of the two families of objects can be defined automatically on the other one -- so that our constructions could be rephrased into the language of trees.
In the end, the relations between free moments and free cumulants, and the underlying combinatorics will be explained through a process similar to (though different from) the one familiar in the theory of dynamical systems under the name ``arborification'', due to J.~Ecalle \cite{Ecalle}, and originating ultimately from Cayley's work on differential operators; we refer to Section~\ref{sect:splitting} for further details on the subject. 

We introduce in the process a new bialgebra (more precisely an unshuffle bialgebra) structure on the tensor algebra over the linear span of non-crossing partitions. The ``splitting process'' to which the title of the article refers explains how generating series of free moments can be lifted to characters on this bialgebra. More precisely, we show that the usual interpretation of the relation between free moments and cumulants in terms of M\"obius calculus on the lattice of non-crossing partitions is governed by a particular (and relatively arbitrary) choice for the lift of the calculus of free moments from characters on tensor algebras to the analog notion in the framework of non-crossing partitions.
This point of view allows in particular a group-theoretic interpretation of free moments and of their relation with free cumulants (Theorems \ref{thm:Gg} and \ref{keyrell}), paving the way hopefully to a renewed understanding of the theory.

\vspace{0.3cm}

The paper is organized as follows. In the next section we introduce both free and classical cumulants and moments. To make this work more accessible, we recall briefly in section \ref{sect:dendalg} the notion of shuffle algebra. Section \ref{sect:bialgNCpart} presents the core result of this work, i.e., the detailed construction of an unshuffle bialgebra of non-crossing partitions. The last section shows how to lift the relation between free moments and cumulants described on the tensor algebra to its analog on the lattice of  non-crossing partitions.

\vspace{0.3cm}

In the following $\mathbb{K}$ is a ground field of characteristic zero. We also assume any $\mathbb{K}$-algebra $A$ to be associative and unital, if not stated otherwise. The unit in $A$ is denoted $1_A$.  

\vspace{0.4cm}
\noindent {\bf{Acknowledgements}}: The first author is supported by a Ram\'on y Cajal research grant from the Spanish government. The second author acknowledges support from the grant ANR-12-BS01-0017, Combinatoire Alg\'ebrique, R\'esurgence, Moules et Applications. Support by the CNRS GDR Renormalisation is also acknowledged.


\section{Classical versus free moments and cumulants}
\label{sect:CM}

In classical probability theory the relation between moments, $m_n$, and cumulants, $c_n$, can be formulated on the level of exponential generating functions ${M}(t) := 1 + \sum_{i > 0} {m}_i \frac{t^i}{i!}$ and $C(t) := \sum_{i > 0} c_i \frac{t^i}{i!}$
$$
	{M}(t) = \exp\big(C(t)\big).
$$
It implies polynomial relations among the coefficients
\begin{equation}
\label{classical1}
	{m}_n = B_n(c_1,\ldots,c_n),
\end{equation}
where $B_n:=B_n(c_1,\ldots,c_n)$ is the Bell polynomial of order $n$. Recall that these polynomials can be defined recursively by
$$
	B_{n+1} = \sum_{k=0}^n {n \choose k} B_{n-k} c_{k+1}, \quad B_0:=1.
$$ 
We denote the Bell polynomials up to order 5
\begin{eqnarray*}
	B_0 &=& 1\\ 
	B_1 &=& c_1, \
	B_2 = c_1^2 + c_2,\
	B_3 = c_1^3 + 3c_1 c_2 + c_3\\
	B_4 &=& c_1^4 + 6c_1^2 c_2 + 3c_2^2 + 4c_1c_3 + c_4\\
	B_5 &=& c_1^5 + 10c_1^2c_3 + 10c_2c_3 + 10c_1^3c_2 + 15c_1c_2^2 + 5c_1c_4 + c_5.
\end{eqnarray*}

Relation (\ref{classical1}) can be formulated using M\"obius calculus on the lattice of set partitions \cite{Beissinger,Rota,Speed}. Recall that a partition $L$ of a (finite) set $[n]:=\{1,\dots,n\}$ consists of a collection of (non-empty) subsets $L=\{L_1,\ldots,L_b\}$ of $[n]$, called blocks, which are mutually disjoint, i.e., $L_i \cap L_j=\emptyset$ for all $i\neq j$, and whose union $\cup_{i=1}^b L_i=[n]$. By $|L|:=b$ the number of blocks of the partition $L$ is denoted, and $|L_i|$ is the number of elements in the $i$th block $L_i$. Given $p,q\in [n]$ we will write that $p \sim_{L} q$ if and only if they belong to the same block. The lattice of set partitions of $[n]$ is denoted by $P_n$. It has a partial order of refinement: $L \leq K$ if $L$ is a finer partition than $K$. The partition $\hat{1}_n = \{L_1\}$ consists of a single block, i.e., $|L_1|=n$, and is the maximum element in $P_n$. The partition $\hat{0}_n=\{L_1,\dots,L_n\}$ has $n$ singleton blocks, and is the minimum partition in $P_n$. 

In the following we denote set partitions pictorially. Consider for instance
$$
		{\nci\ , \atop \!\!\! 1}
		\qquad\ 
		{\ncii\ , \atop {\!\!\! 1\hspace{0.65cm}\; 2} }
		\qquad\ 
		{\nci \quad \nci \quad \nci\ , \atop \!\!\! \! 1\;\; 2\;\; 3 }
$$
the first represents the singleton $\hat{0}_1=\hat{1}_1=\{1\}$ in $P_1$. The second is the single block partition, i.e., the maximal element $\hat{1}_2=\{1,2\} \in P_2$. Then follows the minimal element in $P_3$, i.e., the partition of $[3]$ into singletons, $\hat{0}_3=\{\{1\},\{2\},\{3\}\}$. In $P_3$ we have other partitions
$$
		\nci \quad \ncii\ ,
		\qquad\
		\ncii \quad \nci\ ,
		\qquad
		\nciiia\ ,
		\qquad\
		\nciii\ .
$$
The first, second, and third represent the partitions $\{\{1\},\{2,3\}\}$, $\{\{1,2\},\{3\}\}$, and $\{\{1,3\},\{2\}\}$, respectively. The last one is the maximal element in $P_3$, that is, the single block partition $\hat{1}_3 = \{1,2,3\}$. In $P_4$, $P_5$ and $P_6$ we find for instance
$$
	\nciiic\ ,
	\quad\
	 \nciiib\ ,
	\quad\
	\nciiiia\ ,
	\quad\
	\nciiiiia\ ,
	\quad\
	\nciiiiii
$$ 
which stand for the following partitions, $\{\{1,4\},\{2, 3\}\}$, $\{\{1,3\},\{2, 4\}\}$ and $\{\{1,5\},\{2\}, \{3\}, \{4\}\}$, $\{\{1,5\},\{2\}, \{3,4\}\}$, and $\{\{1,6\},\{2\}, \{3,5\}, \{4\} \}$, respectively.

Moments ${m}_n$ are then given as a convolution product        
\begin{equation}
\label{classical2}
	{m}_n = (c * \zeta)(\hat{1}_n ),
\end{equation}
where $\zeta$ denotes the zeta-function on the lattice $P_n$. Equivalently, $({m}*\mu )(\hat{1}_n) = c_n$, where $\mu$ is the M\"obius function \cite{Rota}. Here both ${m}$ and $c$ are considered to be multiplicative functions on $P_n$, i.e., ${m}(L):={m}_{|L_1|} \cdots {m}_{|L_b|}$, $L \in P_n$.

The relation between free moments $\mathrm{m}_n$ and free cumulants $\mathrm{k}_l$ is rather different. In terms of generating functions it is given by
$$
	F(t) =K(tF(t)),
$$
where $F(t) := 1 + \sum_{i > 0} \mathrm{m}_i t^i$ and $K(t) := 1 + \sum_{i > 0} \mathrm{k}_i t^i$. The first few low order polynomials are
\begin{eqnarray*}
	\mathrm{m}_0 &=& 1 \\
	\mathrm{m}_1 &=& \mathrm{k}_1,\ 
	\mathrm{m}_2 = \mathrm{k}_1^2 + \mathrm{k}_2\\ 
	\mathrm{m}_3 &=& \mathrm{k}_1^3 + 3\mathrm{k}_1 \mathrm{k}_2 + \mathrm{k}_3\\
	\mathrm{m}_4 &=& \mathrm{k}_1^4 + 6\mathrm{k}_1^2 \mathrm{k}_2 
					+ 2\mathrm{k}_2^2 + 4\mathrm{k}_1\mathrm{k}_3 + \mathrm{k}_4.
\end{eqnarray*}
It can be formulated in terms of non-crossing set partitions. Recall that a set partition $L$ of $[n]$ is called non-crossing if for $p_1,p_2,q_1,q_2 \in [n]$ the following does not occur, 
$$
	1\leq p_1<q_1<p_2<q_2\leq n
$$
and 
$$
	p_1\sim_{L}p_2\nsim_{L}q_1\sim_{L}q_2.
$$
The set of non-crossing partitions of $[n]$ will be denoted by $NC_n$. For example, the partition $\{\{1,3\},\{2, 4\}\}$ 
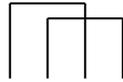
$$
	 \nciiib
$$ 
does not qualify as a non-crossing partition, and therefore is not an element in $NC_4$. One shows that $NC_n$ itself forms a lattice with respect to the partial order $L \leq K$, if $L$ is a finer partition than $K$.  Non-crossing partitions of arbitrary subsets of the integers are defined similarly. We will use repeatedly, and without further notice various elementary properties of non-crossing partitions. In particular, we will use that, if $L$ is a non-crossing partition of $[n]$, then its restriction to an arbitrary subset $S$ of $[n]$ (by intersecting the blocks of $L$ with $S$) defines a non-crossing partition of $S$.

The relation between free moments and cumulants can be formulated using M\"obius calculus on the lattice $NC_n$ 
\begin{equation}
\label{free}
	\mathrm{m}_n = \mathrm{k} * \zeta_{NC},
\end{equation}
where $\zeta_{NC}$ denotes the zeta-function on the lattice $NC_n$, and $\mathrm{m}_n=:\mathrm{m}(\hat{1}_n )$. For more details we refer the reader to \cite{MC,NiSp,Simion}.

\smallskip

For later use we introduce the notion of a hierarchy tree associated to a non-crossing partition. Observe that the blocks of each non-crossing partition are either strictly nested or disjoint. Let us look for example at the non-crossing partition $\{ \{1,4\}, \{2,3\},  \{5,6,7\}\}$
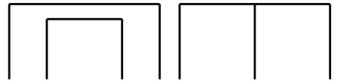
$$
	\nciiic \;\ \nciii\ .
$$
We have that the block $\{2,3\}$ is sitting inside the block $\{1,4\}$. The block $\{5,6,7\}$ is disjoint from the block $\{1,4\}$. To each non-crossing partition we may associate a planar rooted tree $\rho: \mathcal{NC} \to \mathcal{T}_{pl}$, which encodes the hierarchy of the blocks of the non-crossing partition. For the above example we see that $\rho$ maps \scalebox{0.4}{\nciiic \;\ \nciii} to the tree
$$
	\arbref
$$
That is, mark each block with a dot and connect dots by edges according to the nestings of blocks. The deepest sitting blocks are the leaves of the tree. The edges are oriented towards the leaves. The root vertex (green), which has no incoming edge, and is added without being attached to a block, is connected to the outermost disjoint blocks. Note that, of course, this association of trees and non-crossing partitions is highly non-unique, and is not the aforementioned bijection between planar trees and non-crossing partitions. Indeed, the following non-crossing partition 
$$
	\nciiia \;\ \ncii
$$
has the same hierarchy tree as the above one. The following examples of non-crossing partitions
$$
	\ncii \ ,
	\quad\	
	\ncii \;\ \ncii \ ,
	\quad\
	\nciiic\ , 
	\quad\
 	\ncii\;\ \nci \;\ \ncii \ ,
	\quad\
	\nciiic \; \nci\ , 
	\quad\
	\nciiiiia\ , 
	\quad\
$$   
have respectively the following corresponding hierarchy trees
$$
	\arbrea\ ,
	\quad\
	\arbreb\ ,
	\quad\
	\arbrec\ ,
	\quad\
	\arbree\ ,
	\quad\
	\arbref\ ,
	\quad\
	\arbreg\ .
$$ 

\medskip

In \cite{EPnc} the present authors explored a perspective on the relations between (free) cumulants and (free) moments, which differs from the aforementioned approach in terms of M\"obius calculus on $NC_n$. One may argue that the latter approach to free cumulants is not directly related to the physical motivation for their introduction, namely the study of the behaviour of solutions to the random evolution equation following \cite{NeuSpeicher}. We look at the following linear differential equation for the evolution operator $U(t,s)$
$$
	\frac{d U}{d t}(t,t_0)=H(t)U(t,t_0),
$$
where $H(t)$ is a random operator. Writing $< \ \ >$ for the averaging operator, free cumulants are then introduced through a master equation \cite{NeuSpeicher}
\begin{eqnarray*}
	\frac{d}{dt}<U(t,s)> &=& \sum\limits_{i=0}^\infty\;\; \idotsint \limits_{t \geq t_1\geq \dots\geq t_i\geq s}
	\mathrm{k}_{i+1}(t,t_1,\ldots,t_i)<U(t,t_1)> <U(t_1,t_{2})> \cdots <U(t_i,s)> \prod_{j=1}^i dt_j.
\end{eqnarray*}
Using functional calculus this yields the defining recursion:
\begin{eqnarray*}
	<H(t)H(t_1)\cdots H(t_n)> &=&\sum\limits_{r=0}^n \sum\limits_{\{i(1),\ldots ,i(r)\}\subset \{1,\ldots ,n\}} 
	\mathrm{k}_{r+1}(t,t_{i(1)},\dots,t_{i(r)}) <H(t_1)\cdots H(t_{i(1)-1})>\\
	& &\hspace{6cm}\cdots <H(t_{i(r)+1})\cdots H(t_n)>
\end{eqnarray*}
from which the expectation values $<H(t_1)\cdots H(t_n)>$ can be deduced as polynomials in the free cumulants $\mathrm{k}_{l}(t_{1},\dots,t_{l})$, $l=1,\ldots,n$. As it turns out, this recursion is best understood by applying the notion of half-unshuffle calculus on the double bar construction introduced in \cite{EPnc}, that will be recalled below for the sake of completeness.

Second, for group theoretical reasons, the naturalness of using a multiplicative approach to encode the relations between classical or free moments and cumulants, as it is done in the M\"obius inversion approach is not obvious. The stochastic differential equation (SDE) approach, for example, shows that moments should indeed behave multiplicatively and as elements of a group (they are associated to the solution of a SDE), whereas cumulants should behave as elements of a Lie algebra (they behave as generators for the master equation). Although these indications may seem vague to the reader, they can be given a theoretical meaning: the generating series of moments is the solution of a recursive linear equation, and defines a character (i.e., group-like element) on a properly defined Hopf algebra, whereas the generating series of cumulants is the generator of this equation and can be shown to behave as an infinitesimal character (primitive, i.e., a Lie-type element).

The aim of the forthcoming developments in this article, is to show how these ideas combined with a combinatorial approach to free cumulants by means of non-crossing partitions (or equivalent combinatorial objetcs) can be encoded in a proper algebraic framework.


\section{Shuffles}
\label{sect:dendalg}

Recall the abstract definition of a shuffle algebra (also called dendrimorphic algebra)\footnote{We prefer the terminology of shuffles, that is more intuitive, and grounded in the classical calculus of products in algebraic topology, where these notions originated.}. It is a $\mathbb{K}$-vector space $D$ together with two binary compositions $\prec$ and $\succ$ subject to the following three axioms
\begin{eqnarray}
	(a\prec b)\prec c  &=& a\prec(b \prec c + b \succ c)        	\label{A1}\\
  	(a\succ b)\prec c  &=& a\succ(b\prec c)   				\label{A2}\\
   	a\succ(b\succ c)  &=& (a \prec b + a \succ b)\succ c        	\label{A3}.
\end{eqnarray}
These relations yield the associative shuffle product 
\begin{equation}
\label{dendassoc}
	a \shuffle  b := a \prec b + a \succ b
\end{equation}	
on $D$. The two products $\prec,\succ$ are called left and right half-shuffles, respectively. Classical examples of shuffle algebras are provided by chain complexes in algebraic topology, by algebras of operator-valued iterated integrals and by combinatorial Hopf algebras, such as the celebrated Malvenuto--Reutenauer algebra -- equipped with the ``shifted shuffle'' product. Recall the useful notation $a^{\prec n}:=a\prec ( a^{\prec^{n-1}}),\ a^{\prec 0}:=1$ and $\exp^{\prec}(a):=\sum_{n\in \Nb}a^{\prec n}$ (where the latter exponential is an abstract rewriting of the ``time-ordered exponential'' of physicists).

In a commutative shuffle algebra the left and right operations are identified $x \succ y = y \prec x,$ which implies commutativity of (\ref{dendassoc}). Classical examples of commutative shuffle algebras are provided by the calculus of scalar-valued iterated integrals (Chen's calculus). Its abstract formulation is given by the tensor algebra equipped with the recursively defined shuffle product
$$
	v_1\otimes \cdots \otimes v_n  \ \shuffle\ w_1\otimes \cdots \otimes w_m 
			:= v_1\otimes \cdots \otimes v_n\prec w_1\otimes \cdots \otimes w_m
					+ v_1\otimes \cdots \otimes v_n\succ w_1\otimes \cdots \otimes w_m
$$
where
\begin{eqnarray*}
	v_1\otimes \cdots \otimes v_n\prec w_1\otimes \cdots \otimes w_m 
			&:=& v_1\otimes (v_2\otimes \cdots \otimes v_n \ \shuffle\ w_1\otimes  \cdots \otimes w_m),\\
	v_1\otimes \cdots \otimes v_n\succ w_1\otimes \cdots \otimes w_m
			&:=& w_1 \otimes (v_1\otimes \cdots \otimes v_n \ \shuffle\ w_2\otimes \cdots \otimes w_m).
\end{eqnarray*}

The algebraic notion dual to the shuffle product, although familiar in the theory of free Lie algebras \cite{MR,Reutenauer}, has been considered only recently from an abstract, i.e., axiomatic, point of view. We refer the reader to Foissy's seminal work \cite{Foissy}.  

\begin{defn}
A counital unshuffle coalgebra (or counital codendrimorphic coalgebra) is a coaugmented coassociative coalgebra $\overline C = C\oplus \Kb$ with coproduct
\begin{equation}
\label{codend}
	\Delta(c) := \bar\Delta(c) + c \otimes \un + \un \otimes c,
\end{equation}
such that, on $C$, the reduced coproduct $\bar\Delta = \Delta_{\prec} + \Delta_{\succ}$ with 
\begin{eqnarray}
	(\Delta_{\prec} \otimes I) \circ \Delta_{\prec}   &=& (I \otimes \bar\Delta)\circ \Delta_{\prec}        	\label{C1}\\
  	(\Delta_{\succ} \otimes I) \circ \Delta_{\prec}   &=& (I \otimes \Delta_{\prec})\circ \Delta_{\succ} 	\label{C2}\\
   	(\bar\Delta \otimes  I) \circ \Delta_{\succ}         &=& (I \otimes \Delta_{\succ})\circ \Delta_{\succ}   \label{C3}.
\end{eqnarray}
\end{defn}

\noindent The maps $\Delta_{\prec}$ and $\Delta_{\succ}$ are called respectively left and right half-unshuffles. The definition of a (nonunital) unshuffle coalgebra is obtained by removing the unit, that is, $\bar\Delta$ is  acting on $C$, and has a splitting into two half-coproducts, $\Delta_{\prec}$ and $\Delta_{\succ}$, which obey relations (\ref{C1}), (\ref{C2}) and (\ref{C3}).

\begin{defn}
An unshuffle bialgebra is a unital and counital bialgebra $\overline B=B \oplus \Kb$ with product $\cdot_B$ and coproduct $\Delta$, and a counital unshuffle coalgebra $\bar\Delta = \Delta_{\prec} + \Delta_{\succ}$. Moreover, the following compatibility relations hold 
\begin{eqnarray}
	\Delta^+_{\prec}(a \cdot_B b)  &=& \Delta^+_{\prec}(a)  \cdot_B \Delta(b)      	\label{D1}\\
  	\Delta^+_{\succ}(a \cdot_B b)  &=& \Delta^+_{\succ}(a)  \cdot_B \Delta(b),     \label{D2}
\end{eqnarray}
where
\begin{eqnarray}
	\Delta^+_{\prec}(a)  &:=& \Delta_{\prec}(a) + a \otimes \un     	\label{D3}\\
  	\Delta^+_{\succ}(a)  &:=& \Delta_{\succ}(a) + \un \otimes a.     	\label{D4}
\end{eqnarray}
\end{defn}

We refer to classical sources for the definitions of coalgebras, bialgebras and Hopf algebras, see e.g., \cite{Cartier}. Recall that the two notions of bialgebras and Hopf algebras (bialgebras equipped with an antipode) identify when suitable connectedness hypothesis hold (i.e.~when the bialgebra/Hopf algebra can be equipped with a graduation such that the degree zero component is the ground field). This hypothesis will hold in all the examples we will consider, so that there is no difference in this article between the two notions of bialgebras and Hopf algebras.

The classical example of an unshuffle bialgebra is given by the tensor algebra $\overline T(X)$ over an alphabet $X=\{x_1,\ldots ,x_l,\ldots \}$. Notice, for later use, that we write $T(X)$ for the non unital tensor algebra and $\overline T(X)$ for the unital tensor algebra. Both are the linear span of words $x_{i_1}\cdots x_{i_n}$ over $X$ (the empty word is allowed in the unital case, but not in $T(X)$), equipped with the concatenation product of words: $x_{i_1}\cdots x_{i_n} \cdot x_{j_1}\cdots x_{j_m}:=x_{i_1}\cdots x_{i_n}x_{j_1}\cdots x_{j_m}$. The unital tensor algebra $\overline T(X)$ is equipped with the unshuffling coproduct
$$
	\Delta(x_{i_1}\cdots x_{i_n}) = \sum_{I \coprod J = [n]} x_I\otimes x_J,
$$
where, for the subset $S=\{s_1,\ldots ,s_k\}\subset [n]$, $x_S$ stands for the word $x_{i_{s_1}}\cdots x_{i_{s_k}}$. 
Setting 
$$
	\Delta_\prec^+(x_{i_1} \cdots x_{i_n})
		=\sum_{I\coprod J=[n] \atop 1\in I}x_I\otimes x_J
$$ 
defines an unshuffle bialgebra structure on $\overline T(X)$, whose fine structure is studied in \cite{foissypatras}. 

A key example for our purposes is given in terms of the double tensor algebra (or double bar construction) over an associative algebra. Let $A$ be an associative $\mathbb{K}$-algebra. Define $T(A):=\oplus_{n > 0} A^{\otimes n}$ to be the nonunital tensor algebra over $A$. The full tensor algebra is denoted $\overline T(A):=\oplus_{n \ge 0} A^{\otimes n}$. Elements in $T(A)$ are written as words $a_1\cdots a_n \in T(A)$. Note that we denote $a \cdots a\in A^{\otimes n}$ by $a^{\otimes n}$, and the product of the $a_i$s in $A$ is written $a_1 \cdot_A a_2$. Concatenation of words makes $T(A)$ an algebra, which is naturally graded by the length of a word, i.e., its number of letters.

We also set $T(T(A)):=\oplus_{n > 0} T(A)^{\otimes n}$, and use the bar-notation to denote elements $w_1 | \cdots | w_n \in T(T(A))$, $w_i \in T(A)$, $i=1,\ldots,n$. The algebra $T(T(A))$ is equipped with the concatenation product. For $a= w_1 | \cdots | w_n$ and $b=  w_1' | \cdots | w_m'$ we denote their product in $T(T(A))$ by $a|b$, that is, $a|b := w_1 | \cdots | w_n | w_1' | \cdots | w_m'$. This algebra is multigraded, $T(T(A))_{n_1,\ldots ,n_k}:=T_{n_1}(A)\otimes \cdots \otimes T_{n_k}(A)$, as well as graded. The deegre $n$ part is  $T(T(A))_n:=\bigoplus\limits_{n_1+ \cdots +n_k=n}T(T(A))_{n_1,\ldots ,n_k}$. Similar observations hold for the unital case, that is, $\overline T(T(A))=\oplus_{n \ge 0} T(A)^{\otimes n}$, and we will identify without further comments a bar symbol such as $w_1|1|w_2$ with $w_1|w_2$ (formally, using the canonical map from $\overline T(\overline T(A))$ to $\overline T(T(A))$).

\medskip

Given two (canonically ordered) subsets $S \subseteq U$ of the set of integers $\bf N$, we call connected component of $S$ relative to $U$ a maximal sequence $s_1, \ldots , s_n$ in $S$, such that there are no $ 1 \leq i < n$ and $t\in U$, such that $s_i < t <s_{i+1}$. In particular, a connected component of $S$ in $\bf N$ is simply a maximal sequence of successive elements $s,s+1,\ldots ,s+n$ in $S$.

Consider a word $a_1\cdots a_n \in T(A)$. For $S:=\{s_1,\ldots, s_p\} \subseteq [n]$, we set $a_S:= a_{s_1} \cdots a_{s_p}$ (resp. $a_\emptyset:=1$). Denoting $J_1,\ldots,J_k $ the connected components of $[n] - S$, we set $a_{J^S_{[n]}}:= a_{J_1} | \cdots | a_{J_k}$. More generally, for $S \subseteq U \subseteq [n]$, set  $a_{J^S_U}:= a_{J_1} | \cdots | a_{J_k}$, where the $a_{J_j}$ are now the connected components of $U-S$ in $U$.

\begin{defn} \label{def:coproduct}
The map $\Delta : T(A) \to \overline T(A) \otimes  \overline T(T(A))$ is defined by
\begin{equation}
\label{HopfAlg}
	\Delta(a_1\cdots a_n) := \sum_{S \subseteq [n]} a_S \otimes  a_{J_1} | \cdots | a_{J_k}
					   =\sum_{S \subseteq [n]} a_S \otimes a_{J^S_{[n]}}.
\end{equation} 
The coproduct is then extended multiplicatively to all of $\overline T(T(A))$
$$
	\Delta(w_1 | \cdots | w_m) := \Delta(w_1) \cdots \Delta(w_m),
$$
with $\Delta(\un):= \un \otimes \un$.
\end{defn}

The proofs of the following two theorems appeared in \cite{EPnc}.
 
\begin{thm} 
\label{thm:HA}
The graded algebra $\overline T(T(A))$ equipped with the coproduct \eqref{HopfAlg} is a connected graded noncommutative and noncocommutative Hopf algebra. 
\end{thm}

The crucial observation is that coproduct \eqref{HopfAlg} can be split into two parts as follows. On $T(A)$ define the {\it{left half-coproduct}} by
\begin{equation}
\label{HAprec+}
	\Delta^+_{\prec}(a_1 \cdots a_n) := \sum_{1 \in S \subseteq [n]} a_S \otimes a_{J^S_{[n]}},
\end{equation}
and
\begin{equation}
\label{HAprec}
	\Delta_{\prec}(a_1 \cdots a_n) := \Delta^+_{\prec}(a_1 \cdots a_n) - a_1 \cdots a_n \otimes \un. 
\end{equation}
The {\it{right half-coproduct}} is defined by
\begin{equation}
\label{HAsucc+}
	\Delta^+_{\succ}(a_1 \cdots a_n) := \sum_{1 \notin S \subset [n]} a_S \otimes a_{J^S_{[n]}}
\end{equation}
and
\begin{equation}
\label{HAsucc}
	\Delta_{\succ}(a_1 \cdots a_n) := \Delta^+_{\succ}(a_1 \cdots a_n) -  \un \otimes a_1 \cdots a_n.
\end{equation}
Which yields $\Delta = \Delta^+_{\prec} + \Delta^+_{\succ}$, and 
$$
	\Delta(w) = \Delta_{\prec}(w) + \Delta_{\succ}(w) + w \otimes \un + \un \otimes w.
$$
This is extended to $T(T(A))$ by defining
\begin{eqnarray*}
	\Delta^+_{\prec}(w_1 | \cdots | w_m) &:=& \Delta^+_{\prec}(w_1)\Delta(w_2) \cdots \Delta(w_m) \\
	\Delta^+_{\succ}(w_1 | \cdots | w_m) &:=& \Delta^+_{\succ}(w_1)\Delta(w_2) \cdots \Delta(w_m). 
\end{eqnarray*}

\begin{thm}  \cite{EPnc} \label{thm:bialg}
The bialgebra $\overline T(T(A))$ equipped with $\Delta_{\succ}$ and $\Delta_{\prec}$ is an unshuffle bialgebra. 
\end{thm}

Furthermore, recall that the set of linear maps, $Lin(T(T(A)),\mathbb{K})$, is a $\mathbb{K}$-algebra with respect to the convolution product defined for $f,g\in Lin(T(T(A)),\mathbb{K})$ by
$$
	f * g := m_\Kb\circ (f\otimes g)\circ \Delta,
$$ 
where $m_\Kb$ stands for the product map in $\mathbb{K}$. We define accordingly the left and right half-convolution products:
$$
	f\prec g:=m_\mathbb{K}\circ (f\otimes g)\circ \Delta_\prec ,
	\quad\
	{\rm{and}}
	\quad\
	f\succ g:=m_\mathbb{K}\circ (f\otimes g)\circ \Delta_\succ .
$$

\begin{prop}  \cite{EPnc}
The space $({\mathcal L}_A:=Lin(T(T(A)),\mathbb{K}), \prec, \succ)$ is a shuffle algebra.
\end{prop}

The relation between free moments and free cumulants then reads:

\begin{thm}  \cite{EPnc}\label{tim:freeprob}
Let $\phi: A \to \mathbb{K}$ be a unital map, we still write $\phi$ for its extension to $T(A)$ ($\phi(a_1\cdots a_n):=\phi(a_1\cdot_A ... \cdot_A a_n)$). The extension of $\phi$ to $\overline T(T(A))$ is denoted $\Phi$, that is, $\Phi(\omega_1|\cdots |\omega_n):=\phi(\omega_1)\cdot_A ... \cdot_A  \phi(\omega_n)$. Let the map $\kappa: \overline T(T(A))\to \mathbb{K}$ be the solution to $\Phi = e + \kappa\prec \Phi$. For $a \in A$ we set $\mathrm{k}_n:=\kappa(a^{\otimes n})$, $n \geq 1$ and $\mathrm{m}_n:=\Phi(a^{\otimes n})=\phi(a^n)$, $n\geq 0$. Then
$$
	\mathrm{m}_n=\sum\limits_{s=1}^n\sum\limits_{i_1 + \cdots + i_s = n-s} 
					\mathrm{k}_s \mathrm{m}_{i_1} \cdots \mathrm{m}_{i_s}.
$$
In particular, the $\mathrm{k}_n$ identify with the free cumulants of $a \in (A,\phi)$.
\end{thm}

These results generalize to free cumulants in several variables. For $L=\{L_1,\ldots,L_k\} \in NC_n$, $a_1,\ldots,a_n\in A$ and $\omega$ a linear form on $T(A)$, we write $\omega^L(a_1,\ldots ,a_n):= \prod_{i=1}^k\omega(a_{L_i})$. The generalized non-crossing cumulants $R(a_1,\ldots,a_n)$ associated to a unital map $\phi: A \to \mathbb{K}$ are then the multilinear maps defined by the implicit equations that can be solved recursively:
\begin{equation}
\label{stdrel}
	\phi(a_1\cdots a_n)=:\sum\limits_{L\in NC_n}R^L(a_1,\ldots ,a_n).
\end{equation}

\begin{thm}\label{gencumulant}
Let $\phi: A \to \mathbb{K}$ be a unital map, and $\Phi$ its extension to $\overline T(T(A))$  as above. Let the map $\kappa: \overline T(T(A))\to \mathbb{K}$ be the solution to $\Phi = e + \kappa\prec \Phi$. For $a_1,\ldots ,a_n \in A$, we have: $\kappa(a_1 \cdots a_n)=R(a_1,\ldots ,a_n)$. That is, $\kappa$ computes the generalized non-crossing cumulants associated to $\phi$.
\end{thm}

The next sections aim at giving complementary approach to the links between non-crossing partitions and free cumulants -- namely, we will show that these links can be understood through the filter of Hopf algebraic constructions and half-shuffles.

\section{The unshuffle bialgebra of non-crossing partitions}
\label{sect:bialgNCpart}

Let $L=\{L_1,\ldots ,L_k\}$ be an arbitrary non-crossing partition of $[n]:=\{1,\ldots ,n\}$ with $\inf(L_i)<\inf(L_{i+1})$ for $i=1,\ldots ,n-1$. Let us write $L_i < L_j$ if $\forall a \in L_i$ and $\forall b \in L_j$ we have $a<b$. We define a partial order $<_{L}$ on the blocks $L_i$ as follows: $L_i<_{L} L_j$ if and only if, for all $m \in L_i$, $\inf (L_j) < m < \sup(L_j)$. The very definition of non-crossing partitions shows that this partial order is well-defined. Moreover, given two distinct blocks $L_i, L_j\in L$, then one and only one of the following inequalities holds
$$
	L_i < L_j,\quad L_j < L_i,\quad L_i<_L L_j,\quad L_j<_L L_i.
$$

As an example we consider the non-crossing partition $L \in P_{10}$ with 5 block $L=\{L_1,L_2,L_3,L_4,L_5\} =\{ \{1,3,8\},\{2\},\{4,6,7\},\{5\},\{9,10\} \}$
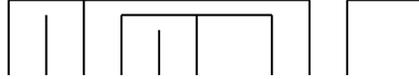
$$
	\nciiiiiia
$$
The block $L_5 > L_i$, $i=1,2,3,4$, and $L_2 <_L L_1$, $L_4 <_L L_3 <_LL_1$.

A partition of the blocks of $L$ into two (possibly empty) subsets
$$
	L = Q \coprod T = \{Q_1,\ldots ,Q_i\} \coprod \{T_1,\ldots ,T_{n-i}\}
$$ 
will be said {\it{admissible}} if and only if for all $p \leq i$, $q \leq n-i$, $Q_p \not<_L T_q$, that is, $T_q <_L Q_p $ or the two subsets of $[n]$ are incomparable for the partial order. We write then $L=Q {\coprod \atop {\rm{adm}}} T$. Admissible partitions of non-crossing partitions of arbitrary finite subsets $S$ of the integers are defined accordingly.
Returning to the above example, we have (the list is not exhaustive)
$$
	L= \{L_1,L_2,L_3,L_4\} {\coprod \atop {\rm{adm}}} \{L_5 \}
	  = \{L_1,L_2,L_5\} {\coprod \atop {\rm{adm}}} \{L_3,L_4\}
	  = \{L_1,L_5\} {\coprod \atop {\rm{adm}}} \{L_2, L_3,L_4\}
$$

Similarly, a partition 
$$
	L = Q \coprod T\coprod U
	    = \{Q_1,\ldots ,Q_i\}\coprod \{T_1,\ldots ,T_{p-i}\}\coprod \{U_1,\ldots ,U_{n-p}\}
$$
of $L$ is admissible if and only if for all $v \leq i,\ w\leq p-i,\ z\leq n-p$, $Q_v \not<_L T_w,\ Q_v\not<_L U_z,\ T_w\not<_L U_z$. Notice that there is a bijection between admissible partitions $L=Q\coprod T\coprod U$, pairs of admissible decompositions $L=Q \coprod W,\ W=T\coprod U$ and pairs of admissible decompositions $L=V\coprod U,\ V=Q\coprod T$. We will refer to this property as the {\it{coassociativity of admissibility}} and write $L=Q\!{\coprod\atop {\rm{adm}}}\!T\!{\coprod\atop {\rm{adm}}}\!U$.

For an admissible partition $L=Q {\coprod \atop {\rm{adm}}} U$ as above, we consider as in the previous section the connected components $J_1,\ldots , J_{k(L ,Q)}$ of $[n]-(Q_1\cup \cdots \cup Q_i)$, that we will call slightly abusively from now on the {\it{connected components}} of $[n]-Q$ . The definition of $<_L$ implies that $J_i\cap U_j$ is empty or equals $U_j$. We write $J_i^{L,Q}$ for the set of all non-empty intersections $J_i \cap U_j,\ j=1, \ldots , n-i$ and notice that, since $L$ is a non-crossing partition of $[n]$, $J_i^{L,Q}$ is, by restriction, a non-crossing partition of the component $J_i$. For the same reason, $Q$ is a non-crossing partition of $Q_1\cup \cdots \cup Q_i$.

Let us recall now that, given a finite subset $S$ of cardinality $n$ of the integers, the standardization map $st$ is the (necessarily unique) increasing bijection between $S$ and $[n]$. By extension, we write also $st$ for the induced map on the various objects associated to $[n]$ (such as partitions). For example, the standardization of the non-crossing partition $L:=\{\{3,6,10\},\{4,5\},\{8\}\}$ of the set $\{3,4,5,6,8,10\}$ is the non-crossing partition $st(L):=\{\{1,4,6\},\{2,3\},\{5\}\}$ of $[6]=st(\{3,4,5,6,8,10\})$.

The linear span $\mathcal{NC}$ of all non-crossing partitions can then be equipped with a coproduct map $\Delta$ from $\mathcal{NC}$ to $\mathcal{NC}\otimes T(\mathcal{NC})$ defined by (using our previous notations and the bar $|$ notation for elements in $T(\mathcal{NC})$)
$$
	\Delta(L)=\sum\limits_{Q\!\coprod\limits_{\rm{adm}}\!U = L} 
	st(Q) \otimes \big(st(J_1^{L,Q})| \cdots | st(J_{k(L, Q)}^{L,Q})\big).
$$
A few examples may be helpful at this stage. 
$$
	\Delta(\{\{1,4\},\{2, 3\}\}) = \{\{1,4\},\{2, 3\}\} \otimes \un + \un \otimes \{\{1,4\},\{2, 3\}\} 
						+ \{1,2\} \otimes \{1,2\}
$$
\begin{eqnarray*}
	\Delta(\{\{1,5\},\{2\}, \{3,4\}\}) &=&\{\{1,5\},\{2\}, \{3,4\}\} \otimes \un + \un \otimes \{\{1,5\},\{2\}, \{3,4\}\}
						+  \{\{1,3\},\{2\}\} \otimes \{1,2\} \\
						& & + \{\{1,4\},\{2, 3\}\} \otimes \{1\}
						+ \{1,2\} \otimes \{\{1\}, \{2,3\}\} 
\end{eqnarray*}
\begin{eqnarray*}
	\Delta(\{\{1,2\}, \{3\},\{4\}\}) &=& \{\{1,2\}, \{3\},\{4\}\} \otimes \un + \un \otimes \{\{1,2\}, \{3\},\{4\}\}
						+ \{1,2\} \otimes \{\{1\}\{2\}\}\\
					       & & + 2 \{\{1,2\},\{3\}\} \otimes \{1\} 
					       		+ \{1\}  \otimes \{\{1,2\},\{3\}\}
							+ \{1\}  \otimes \{1,2\} | \{1\}\\
						& & 	 \{\{1\},\{2\}\}  \otimes \{1,2\}
\end{eqnarray*}
The graphical notation of the coproduct simplifies since it is automatically standardized (the bar of the bar notation is written in bold to distinguish it from the trivial tree)
$$
	\Delta(\scalebox{0.5}{\nci}\ ) = \scalebox{0.5}{\nci} \otimes \un + \un \otimes \scalebox{0.5}{\nci} 
	\qquad\
	\Delta(\scalebox{0.5}{\ncii}\ ) = \scalebox{0.5}{\ncii} \otimes \un + \un \otimes \scalebox{0.5}{\ncii}
$$
$$
	\Delta(\scalebox{0.5}{\nciiia}\ ) = \scalebox{0.5}{\nciiia} \otimes \un + \un \otimes \scalebox{0.5}{\nciiia} 
							\ +\  \scalebox{0.5}{\ncii} \otimes \scalebox{0.5}{\nci}
$$
$$
	\Delta(\scalebox{0.5}{\nciiic}\ ) =  \scalebox{0.5}{\nciiic} \otimes \un + \un \otimes \scalebox{0.5}{\nciiic} 
							\ +\  \scalebox{0.5}{\ncii} \otimes \scalebox{0.5}{\ncii}
$$
\begin{eqnarray*}
	\Delta(\scalebox{0.5}{\nciiiiia}\ ) &=& \scalebox{0.5}{\nciiiiia} \otimes \un + \un \otimes \scalebox{0.5}{\nciiiiia}
								\ + \ \scalebox{0.5}{\nciiia} \otimes \scalebox{0.5}{\ncii} 
								\ + \ \scalebox{0.5}{\nciiic} \otimes \scalebox{0.5}{\nci}
								\ +\ \scalebox{0.5}{\ncii} \otimes \scalebox{0.5}{\nci \ncii}
\end{eqnarray*}
\begin{eqnarray*}
	\Delta(\scalebox{0.5}{\ncii \nci \nci}\ ) &=& \scalebox{0.5}{\ncii \nci \nci} \otimes \un 
									+ \un \otimes \scalebox{0.5}{\ncii \nci \nci}
								\ + \ \scalebox{0.5}{\ncii} \otimes \scalebox{0.5}{\nci \nci} 
								\ + \ \scalebox{0.5}{\nci} \otimes \scalebox{0.5}{\ncii \nci}
								\ + \ \scalebox{0.5}{\nci} \otimes \scalebox{0.5}{\ncii} \;\; \baarr\ \scalebox{0.5}{\nci}
								\ +\ 2 \scalebox{0.5}{\ncii \nci} \otimes \scalebox{0.5}{\nci}
								\ +\  \scalebox{0.5}{\nci \nci} \otimes \scalebox{0.5}{\ncii}
\end{eqnarray*}

This coproduct can be understood on hierarchy trees associated to non-crossing partitions in terms of admissible edge cuts. Recall that each oriented edge has an arrival vertex. We write $E(t)$ for the set of edges of the tree $t$. First, we introduce an elementary edge cut, that is,  the removal from the rooted tree of a single specified edge $e \in E(t)$. The resulting pair of trees consists of the rooted part, $R_{e}$ i.e., the tree with the original root of $t$, and the tree $t^{(1)}$, which has the arrival vertex of the cut edge as root. The latter becomes the pruned tree, $P_{e}$, corresponding to the cut at $e$ by connecting via an edge a new root vertex to the root of the tree $t^{(1)}$. 

Now, an admissible edge cut $c^t$ of a tree $t$ is a collection of edges, $c^t \subset E(t)$, such that any path from the root to a leaf contains at most one edge of the collection. The corresponding rooted and pruned parts are given as follows. 

First, we assume that the admissible edge cut $c^t$ consists of two edges of $t$. These two edges may be outgoing edges of either two different vertices or of the same vertex. In the latter case they may be adjacent or not. Again, from deleting the edges of $c^t$ in $t$, we obtain the rooted part, $R_{c^t}$, and we obtain two trees $t^{(1)}, t^{(2)}$. If the edges were outgoing from different vertices, or if they were non-adjacent at the same vertex, then the resulting pruned part, $P_{c^t}$, consists of a monomial of two rooted trees resulting from adding a new root to $t^{(1)}$ and a new root to $t^{(2)}$. Since $t$ is planar, the order of the monomial follows from the order of the trees $t^{(1)}, t^{(2)}$ in $t$. In the case of the two edges in $c^t$ being adjacent, the resulting pruned part, $P_{c^t}$, consist of a single tree following from adding a single new root to the roots of $t^{(1)}, t^{(2)}$. As an example we look at the coproduct of 
the non-crossing partition \scalebox{0.4}{\ncii \;\ \ncii \;\ \ncii}
$$
	\Delta(\scalebox{0.4}{\ncii \;\ \ncii \;\ \ncii}\ ) 
	= \scalebox{0.4}{\ncii \;\ \ncii \;\ \ncii} \otimes \un + \un \otimes \scalebox{0.4}{\ncii \;\ \ncii \;\ \ncii}
	+ 3 \scalebox{0.4}{\ncii \;\ \ncii} \otimes  \scalebox{0.4}{\ncii}
	+ 2  \scalebox{0.4}{\ncii} \otimes  \scalebox{0.4}{\ncii \;\ \ncii}
	+  \scalebox{0.4}{\ncii} \otimes  \scalebox{0.4}{\ncii} \scalebox{0.8}{\baarr} \, \scalebox{0.4}{\ncii}.
$$
For the corresponding hierarchy tree
$$
	\arbree
$$           
we find that $\tilde{\Delta}\circ \rho(\scalebox{0.4}{\ncii \;\ \ncii \;\ \ncii}\ )$ equals
$$
	\tilde{\Delta}(\scalebox{0.5}{\arbree}) = \scalebox{0.5}{\arbree} \otimes \un + \un \otimes \scalebox{0.5}{\arbree}
	+ 3 \scalebox{0.5}{\arbreb} \otimes  \scalebox{0.5}{\arbrea}
	+ 2  \scalebox{0.5}{\arbrea} \otimes  \scalebox{0.5}{\arbreb}
	+  \scalebox{0.5}{\arbrea} \otimes  \scalebox{0.5}{\arbrea} \scalebox{0.5}{\arbrea},
$$
where we write $\tilde{\Delta}$ for the coproduct induced by $\Delta$ on trees.
The term with the coefficient two corresponds to the two admissible cuts at two adjacent edges. The last term corresponds to the admissible cut with two cut edges, which are not adjacent. 
Another example is the coproduct of the tree
$$
	\tilde\Delta(\scalebox{0.5}{\arbref}) = \scalebox{0.5}{\arbref} \otimes \un + \un \otimes \scalebox{0.5}{\arbref}
	+ \scalebox{0.5}{\arbrec} \otimes \scalebox{0.5}{\arbrea}
	+ \scalebox{0.5}{\arbreb} \otimes \scalebox{0.5}{\arbrea}
	+ \scalebox{0.5}{\arbrea} \otimes \scalebox{0.5}{\arbrec}
	+ \scalebox{0.5}{\arbrea} \otimes \scalebox{0.5}{\arbrea} \scalebox{0.5}{\arbrea}
$$ 
which would correspond for instance to the coproduct of the non-crossing partition \scalebox{0.4}{\nciiia\ \ncii}.

For a general admissible cut $c^t$ we have that adjacent cut edges result in a single rooted tree in the pruned part $P_{c^t}$, non-adjacent cut edges result in single rooted trees in the monomial $P_{c^t}$. 

\smallskip

The map $\Delta$ is then extended multiplicatively to a coproduct on $\overline T(\mathcal{NC})$ 
$$
	\Delta (L_1|\cdots |L_n):=\Delta(L_1)\cdots \Delta(L_n), \ \Delta(1)=1\otimes 1.
$$
$\overline T(\mathcal{NC})$ is equipped with the structure of a free associative algebra over $\mathcal{NC}$ by the concatenation map, $(L_1|\cdots |L_k)\cdot (L_{k+1}|\cdots |L_n):=(L_{1}|\cdots |L_k | L_{k+1}| \cdots |L_n)$.

\begin{thm} 
\label{thm:HANC}
The graded algebra $\overline T(\mathcal{NC})$ equipped with the coproduct $\Delta$ is a connected graded noncommutative and noncocommutative Hopf algebra. 
\end{thm}

Here, the grading is the obvious one, i.e., a non-crossing partition of $[n]$ is considered to be of degree $n$. Since $\Delta$ is, by definition, a multiplicative map from $\overline T(\mathcal{NC})$ to $\overline T(\mathcal{NC})\otimes \overline T(\mathcal{NC})$, proving the Theorem amounts to proving coassociativity of $\Delta$. 
We have, writing $\id$ for the identity map:
\begin{eqnarray*}
	     A &:=&(\Delta \otimes \id) \circ \Delta (L) = (\Delta\otimes \id) \big( \sum\limits_{Q\!\coprod\limits_{\rm{adm}}\!W=L}
	     		st(Q)\otimes (st(J_1^{L,Q})| \cdots |st(J_{k(L, Q)}^{L,Q}))\big)\\
		&=&\sum\limits_{Q\!\coprod\limits_{\rm{adm}}\!W =L}\ \sum\limits_{ U \!\coprod\limits_{\rm{adm}}\!V=Q}
			st(U)\otimes \big(st(J_1^{Q,U})|\cdots |st(J_{k(Q,U)}^{Q,U})\big)\otimes \big(st(J_1^{L,Q})|\cdots |st(J_{k(L, Q)}^{L,Q})\big)
\end{eqnarray*}
where we used that the admissible decompositions of $Q$ are in bijection with the admissible decompositions of $st(Q)$, and where $J_1^{Q,U},\ldots ,J_{k(Q,U)}^{Q,U}$ are the non-crossing partitions associated to the decomposition $U\!\coprod\limits_{\rm{adm}}\!V=Q$.

From the coassociativity of admissibility, we get
\begin{eqnarray*}
	A&=&\sum\limits_{U\!\coprod\limits_{\rm{\rm{adm}}}\! V\!\coprod\limits_{\rm{\rm{adm}}}\!W =L}
		st(U)\otimes \big(st(J_1^{\mbox{\tiny{$ U\!\!\coprod\limits_{adm}\!\!V,U$}}})| \cdots |st(J_{k({\mbox{\tiny{$ U\!\!\coprod\limits_{adm}\!\!V,U$}}}}^{{\mbox{\tiny{$ U\!\!\coprod\limits_{adm}\!\!V,U$}}}})\big)\\
	& & \hspace{4cm} \otimes \big(st(J_1^{{\mbox{\tiny{$L,U\!\coprod\limits_{\rm{adm}}\!V$}}}})| \cdots |st(J_{k({\mbox{\tiny{$L,U\!\coprod\limits_{\rm{adm}}\!V$}}})}^{{\mbox{\tiny{$L,U\!\coprod\limits_{\rm{adm}}\!V$}}}})\big).
\end{eqnarray*}
On the other hand, 
$$
	B:=(\id\otimes \Delta)\circ\Delta (L) = (\id\otimes \Delta)\big(\sum\limits_{U\!\coprod\limits_{\rm{adm}}\!W=L}
		st(U)\otimes (st(J_1^{L,U})| \cdots |st(J_{k(L,U)}^{L,U}))\big).
$$
However, since $J_1^{L,U}\cup \cdots \cup J_{k(L,U)}^{L,U}=W$, with $J_i^{L,U}<J_{i+1}^{L,U}$, families of admissible decompositions of the $J_i^{L,U}$ are in bijection with the admissible decompositions of $W$. Taking into account that standardization commutes with admissible decompositions, we get (with a self-explaining notation for $J_i^{L,U}\cap V$)
\begin{eqnarray*}
	& &\Delta(st(J_1^{L,U}))\cdot \cdots \cdot\Delta(st(J_{k(L,U)}^{L,U}))
		=\sum\limits_{V\!\coprod\limits_{\rm{adm}}\!E=W} 
		\big(st(J_1^{L,U}\cap V)|\cdots | st(J_{k(L,U)}^{L,U}\cap V)\big)\otimes \\	
	& &\big[
		(st(J_1^{J_1^{L,U},J_1^{L,U}\cap V})|\cdots |st(J_{k(J_1^{L,U},J_1^{L,U}\cap V)}^{J_1^{L,U},J_1^{L,U}\cap V}))|\cdots 
		| (st(J_1^{J_{k(L,U)}^{L,U},J_{k(L,U)}^{L,U}\cap V})|\cdots 
		|st(J_{k(J_{k(L,U)}^{L,U},J_{k(L,U)}^{L,U}\cap V)}^{J_{k(L,U)}^{L,U},J_{k(L,U)}^{L,U}\cap V}))\big].
\end{eqnarray*}
In this formula, $J_1^{L,U}\cap V$ is the non-crossing partition in the intersection of the first connected component of $[n]-U$ with $V$, and identifies therefore with the non-crossing partition in the first connected component of $V$ in $U\cup V$, that is, $J_1^{{\mbox{\tiny{$U\!\coprod\limits_{\rm{adm}}\!V,U$}}}}$, and similarly for the other $J_i^{L,U}\cap V$. Similarly, $J_1^{J_1^{L,U},J_1^{L,U}\cap V}$ is the non-crossing partition in the first connected component of $[n] -(U\coprod V)$ and identifies therefore with $J_1^{{\mbox{\tiny{$L,U\!\coprod\limits_{\rm{adm}}\!V$}}}}$, and similarly for the other components in the leftmost right hand side of the last expansion of $B$.
Using again the coassociativity of admissibility, we conclude
\begin{eqnarray*}
	B=\sum\limits_{U\!\coprod\limits_{\rm{adm}}\! V\!\coprod\limits_{\rm{adm}}\!W}
	st(U) \otimes (st(J_1^{{\mbox{\tiny{$U\!\coprod\limits_{\rm{adm}}\!V,U$}}}})| \cdots 
	|st(J_{k({\mbox{\tiny{$U\!\coprod\limits_{\rm{adm}}\!V,U$}}})}^{{\mbox{\tiny{$U\!\coprod\limits_{\rm{adm}}\!V,U$}}}}))\\
	\otimes (st(J_1^{{\mbox{\tiny{$L,U\!\coprod\limits_{\rm{adm}}\!V$}}}})| \cdots |st(J_{k({\mbox{\tiny{$L,U\!\coprod\limits_{\rm{adm}}\!V$}}})}^{{\mbox{\tiny{$L,U\!\coprod\limits_{\rm{adm}}\!V$}}}})),
\end{eqnarray*}
so that $A=B$ and the claim of the theorem follows.

\smallskip

This coproduct can be split into two parts as follows. On $\mathcal{NC}$ define the {\it{left half-coproduct}} by 
\begin{equation}
\label{HANCprec+}
	\Delta^+_{\prec}(L)=\sum\limits_{Q\!\coprod\limits_{\rm{adm}}\!U=L \atop 1 \in Q_1}
					st(Q)\otimes \big(st(J_1^{L,Q})| \cdots |st(J_{k(L , Q)}^{L,Q})\big),
\end{equation}
and
\begin{equation}
\label{HANCprec}
	\Delta_{\prec}(L) := \Delta^+_{\prec}(K) - L \otimes \un. 
\end{equation}
The {\it{right half-coproduct}} is defined by
\begin{equation}
\label{HANCsucc+}
	\Delta^+_{\succ}(L)=\sum\limits_{Q\!\coprod\limits_{\rm{adm}}\!U=L \atop 1\in T_1}
					st(Q)\otimes \big(st(J_1^{L,Q})| \cdots |st(J_{k(L , Q)}^{L,Q})\big),
\end{equation}
and
\begin{equation}
\label{HANCsucc}
	\Delta_{\succ}(L) := \Delta^+_{\succ}(L) -  \un \otimes L .
\end{equation}
Which yields $\Delta = \Delta^+_{\prec} + \Delta^+_{\succ}$, and for $L \in \mathcal{NC}$
$$
	\Delta(L) = \Delta_{\prec}(L) + \Delta_{\succ}(L) + L \otimes \un + \un \otimes  L.
$$
This is extended to $\overline T(\mathcal{NC})$ by defining
\begin{eqnarray}
\label{extens}
	\Delta^+_{\prec}(L_1 | \cdots | L_m) &:=& \Delta^+_{\prec}(L_1)\Delta(L_2) \cdots \Delta(L_m) \\\label{extens2}
	\Delta^+_{\succ}(L_1 | \cdots | L_m) &:=& \Delta^+_{\succ}(L_1)\Delta(\pi_2) \cdots \Delta(L_m). 
\end{eqnarray}

\begin{thm} \label{thm:bialgNC}
The bialgebra $\overline T(\mathcal{NC})$ equipped with $\Delta_{\succ}$ and $\Delta_{\prec}$ is an unshuffle bialgebra. 
\end{thm}

The proof is very similar to the one of the coassociativity of $\Delta$, and details are therefore left to the reader. We will just sketch it. 

The two equations expressing $\Delta_\prec^+(a | b)$ and $\Delta_\succ^+(a | b)$ in the definition of unshuffle bialgebras are automatically satisfied, due to equations (\ref{extens}) and (\ref{extens2}). It remains to show that $\Delta_\prec$ and $\Delta_\succ$ satisfy the axioms of an unshuffle coalgebra. Let us indicate for example how the identity $(\Delta_\prec\otimes I)\circ \Delta_\prec =(I\otimes\overline\Delta)\circ \Delta_\prec$ follows. Following the same arguments as in the proof of coassociativity of $\Delta$, we get
\begin{eqnarray*}
	(\Delta_\prec\otimes I)\circ \Delta_\prec(L)&=&
\sum\limits_{U\!\coprod\limits_{\rm{adm}}\! V\!\coprod\limits_{\rm{adm}}\!W=L \atop 1\in U, V\not=\emptyset, W \not=\emptyset}
	st(U) \otimes (st(J_1^{{\mbox{\tiny{$U\!\coprod\limits_{\rm{adm}}\!V,U$}}}})| \cdots |st(J_{k({\mbox{\tiny{$U\!\coprod\limits_{\rm{adm}}\!V,U$}}})}^{{\mbox{\tiny{$U\!\coprod\limits_{\rm{adm}}\!V,U$}}}})) \otimes\\ 
		& &(st(J_1^{{\mbox{\tiny{$L,U\!\coprod\limits_{\rm{adm}}\!V$}}}})| \cdots |st(J_{k({\mbox{\tiny{$L,U\!\coprod\limits_{\rm{adm}}\!V$}}})}^{{\mbox{\tiny{$L,U\!\coprod\limits_{\rm{adm}}\!V$}}}}))
		=(\id \otimes\overline\Delta)\circ \Delta_\prec(L).
\end{eqnarray*}
The other identities follow from similar computations.

\section{The splitting process}
\label{sect:splitting}

The present section aims at explaining, from a Hopf-theoretic perspective, the transition from the ``double bar construction'' point of view on free moments and cumulants, as introduced earlier in this article, to the non-crossing partitions one.

As alluded at in the introduction, this process is similar in various respects to Ecalle's arborification process -- a technique particularly well suited for tackling problems related to diophantian small denominators in the theory of dynamical systems \cite{Ecalle}. This is due partly to the fact that non-crossing partitions are in bijection with planar rooted trees, but there is more to it. Recall that the very idea of arborification (and of the dual coarborification process) \cite{Menous} is, roughly, to replace computations (e.g. of conjugating diffeomorphisms in the study of normal forms of differential equations) involving formal power series in noncommuting variables (the essence of Ecalle's ``mould calculus'': the noncommuting variables are naturally associated to differential operators \cite{sauzin,menous2}) by computations involving formal power series parametrized by decorated trees. As shown in \cite{Menous}, this process is best understood using the point of view of characters on shuffle and quasi-shuffle algebras.

Although the combinatorial Hopf algebras we deal with are different and carry more structure (we take advantage of the existence of half-shuffles and unshuffles, whereas arborification technique do not), the very idea of the splitting process to be introduced now is formally similar to the one underlying arborification.

For $A$ an arbitrary associative algebra, let us write ${\mathcal{NC}}(A)$, and call the set of $A$-decorated (or simply decorated) non-crossing partitions the graded vector space $\Kb\oplus \bigoplus_{n\in\Nb_\ast}{\mathcal{NC}_n}\otimes(A^{\otimes n})$. For $S=\{s_1,\ldots,s_k\}\subset[n]$ and $a:=a_1\cdots a_n$ (we use the word notation for the elements of $A^{\otimes n}$), let us write once again $a_S:=a_{s_1} \cdots  a_{s_k}$. Similarly, for $Q=\{Q_1,\ldots,Q_i\}$ a non-crossing partition of $S$, we write $a_Q:=a_S$.

The definitions in the previous paragraph carry over to decorated non-crossing partitions in a straightforward way. 
For example, the coproduct map $\Delta$ is defined on ${\mathcal{NC}}(A)$ by
\begin{eqnarray*}
	\Delta(L\otimes (a_1 \cdots a_n))=
	\sum\limits_{Q\!\coprod\limits_{\rm{adm}}\!U=L}
	\big(st(Q)\otimes a_Q\big) \otimes \big(st(J_1^{L,Q})\otimes (a_{J_1^{L,Q}})| \cdots |st(J_{k(L, Q)}^{L,Q})
	\otimes (a_{J_{k(L, Q)}^{L,Q}})\big),
\end{eqnarray*}
where $L$ is a non-crossing partition of $[n]$. It is then extended to  $T({\mathcal{NC}}(A))$ multiplicatively as in the previous sections; the other structural maps on  $T({\mathcal{NC}})$ are extended similarly to $T({\mathcal{NC}}(A))$ and are written with the same symbols are previously. We obtain finally

\begin{thm} \label{thm:bialgNCA}
The bialgebra $\overline T(\mathcal{NC}(A))$ equipped with $\Delta_{\succ}$ and $\Delta_{\prec}$ is an unshuffle bialgebra. 
\end{thm}

\begin{defn}
The splitting map $Sp$ is the map from $T(A)$ to $\mathcal{NC}(A)$ defined by:
$$
	Sp(a_1 \cdots a_n):=\sum\limits_{L \in {NC}_n}L\otimes (a_1 \cdots a_n).
$$
It is extended multiplicatively to a unital map $Sp$ from $\overline T(T(A))$ to $\overline T(\mathcal{NC}(A))$: for $A_1,\ldots ,A_k\in T(A)$,
$$
	Sp(A_1|\cdots |A_k):=(Sp(A_1)|\cdots |Sp(A_k)).
$$
\end{defn} 

The name ``splitting'' map is chosen because, as we shall see later, on dual spaces it permits to ``split'' the value of a linear form on $T(A)$, $\phi(a_1\cdots a_n)$, into a sum indexed by non-crossing partitions.

\begin{thm} \label{thm:splitHom}
The map $Sp$ from $\overline T(T(A))$ to $\overline T(\mathcal{NC}(A))$ is an unshuffle bialgebra morphism.
\end{thm}

The map is multiplicative by its very definition. Let us show for example that it commutes with the coproduct map (the commutation with other structure maps follows by the same arguments).

We have from (\ref{def:coproduct}) that:
$$
	\Delta(a_1\cdots a_n) := \sum_{S \subseteq [n]} a_S \otimes  a_{J_1} | \cdots | a_{J_k}.
$$					

On the other hand: 
\begin{eqnarray*}
	A&:=&\Delta\big(Sp(a_1\cdots a_n)\big)
		=\Delta(\sum\limits_{L \in {NC}_n} L \otimes (a_1\otimes \cdots \otimes a_n))\\
		&=& \sum\limits_{L \in {NC}_n}\sum\limits_{Q\!\coprod\limits_{\rm{adm}}\!U=L}
		   \big(st(Q)\otimes a_Q\big)\otimes \big(st(J_1^{L,Q})\otimes (a_{J_1^{L,Q}})|\cdots |st(J_{k(L , Q)}^{L,Q})\otimes (a_{J_{k(L , Q)}^{L,Q}})\big).
\end{eqnarray*}

Let us consider in the last expression the partial sum corresponding to the set of all $Q$ that are non-crossing partitions of $S\subseteq [n]$. The $J_{i}^{L,Q}$ are then non-crossing partitions of the connected components $J_i$ of $[n]-S$. There are no other restrictions on $Q$ and the $J_{i}^{L,Q}$ and we get:
\begin{eqnarray*}
	A&=&\sum\limits_{S\subseteq [n]}\sum\limits_{Q\in NC(S)}\sum\limits_{J_i^{L,Q}\in NC(J_i)} (st(Q)\otimes a_Q)\otimes \\
	  & & \big(st(J_1^{L,Q})\otimes (a_{J_1^{L,Q}})| \cdots |st(J_{k(L, Q)}^{L,Q})\otimes (a_{J_{k(L, Q)}^{L,Q}})\big),
\end{eqnarray*}
from which the identity with 
$\sum_{S \subseteq [n]} Sp(a_S) \otimes  Sp(a_{J_1}) | \cdots |Sp( a_{J_k})$ follows.

Recall that, since $\overline T(\mathcal{NC}(A))$ is a Hopf algebra, the set of linear maps, $Lin(T(\mathcal{NC}(A)),\Kb)$, is a $\mathbb{K}$-algebra with respect to the convolution product defined in terms of the coproduct $\Delta$, i.e., for $f,g\in Lin(T(\mathcal{NC}(A)),\Kb)$
$$ 
	f \ast g := m_\mathbb{K}\circ (f\otimes g)\circ \Delta,
$$ 
where $m_\mathbb{K}$ stands for the product map in $\mathbb{K}$. Notice that, motivated by the next proposition, we will also use later a shuffle notation for this product (so that $f\ast g=f\shuffle g$). We define accordingly the left and right half-convolution products:
$$
	f\prec g:=m_\mathbb{K}\circ (f\otimes g)\circ \Delta_\prec ,
$$
$$
	f\succ g:=m_\mathbb{K}\circ (f\otimes g)\circ \Delta_\succ .
$$

\begin{prop}
The space $({\mathcal L_{NC}}(A):=Lin(T(\mathcal{NC}(A)),\mathbb{K}), \prec, \succ)$ is a shuffle algebra.
\end{prop}

For completeness, and in view of the importance of this proposition for forthcoming developments, we recall briefly its proof: for arbitrary $f,g,h\in {\mathcal L}_A$,
$$
	(f\prec g)\prec h	=m_\mathbb{K}\circ ((f\prec g)\otimes h)\circ\Delta_\prec 
				=m_\mathbb{K}^{[3]}\circ (f\otimes g\otimes h)\circ(\Delta_\prec\otimes I)\circ \Delta_\prec,
$$
where $m_\mathbb{K}^{[3]}$ stands for the product map from $\mathbb{K}^{\otimes 3}$ to $\mathbb{K}$. Similarly 
\begin{eqnarray*}
	f\prec (g\shuffle h)&=&m_\mathbb{K}\circ (f\otimes (g\shuffle h))\circ \Delta_\prec \\
				&=&m_\mathbb{K}^{[3]}\circ (f\otimes g\otimes h)\circ (I\otimes \overline\Delta)\circ\Delta_\prec,
\end{eqnarray*}
so that the identity $(f\prec g)\prec h=f\prec (g\shuffle h)$ follows from $(\Delta_\prec\otimes I)\otimes\Delta_\prec =(I\otimes \overline\Delta)\circ\Delta_\prec$, and similarly for the other identities characterizing shuffle algebras.

As usual, we equip the shuffle algebras ${\mathcal L_{NC}}(A)$ (and similarly ${\mathcal L}_A:=Lin(T(T(A)),\mathbb{K})$) with a unit. That is, we set $\overline{{\mathcal L_{NC}}(A)}:={\mathcal L_{NC}}(A)\oplus \mathbb{K}\cong Lin(\overline T(\mathcal{NC}(A)),\mathbb{K})$, where in the last isomorphism the unit $\un \in \overline{{\mathcal L_{NC}}(A)}$ is identified with the augmentation map $e \in Lin(\overline T(\mathcal{NC}(A)),\mathbb{K})$ -- the null map on $T(\mathcal{NC}(A))$ and the identity map on $\mathcal{NC}(A)^{\otimes 0}\cong \mathbb{K}$. Moreover, for an arbitrary $f$ in ${\mathcal L_{NC}}(A)$,
$$
	f\prec e=f=e\succ f,\ \ e\prec f=0=f\succ e.
$$

The previous constructions are functorial: the linear dual of an unshuffle bialgebra is a shuffle algebra, and a morphism $f$ between two unshuffle bialgebras induces a morphism of shuffle algebras written $f^\ast$ between the linear duals. 
In particular:

\begin{thm} \label{thm:splitHomdual}
The map $Sp$ from $T(T(A))$ to $T(\mathcal{NC}(A))$ induces a morphism $Sp^\ast$ of shuffle bialgebras with units from  $\overline{{\mathcal L_{NC}}(A)}$ to $\overline{{\mathcal L}_A}$.
\end{thm}

The constructions in \cite{EPnc} of $\overline{{\mathcal L}_A}$ apply similarly to $\overline{{\mathcal L_{NC}}(A)}$.
In particular, let $\phi$ be a linear form on $\mathcal{NC}(A)$. It extends uniquely to a multiplicative linear form $\Phi$ on $T(\mathcal{NC}(A))$ by setting
$$
	\Phi(w_1| \cdots |w_n):=\phi(w_1) \cdots \phi(w_n),
$$
(or to a unital and multiplicative linear form on $\overline T(\mathcal{NC}(A))$). Conversely any such multiplicative map $\Phi$ gives rise to a linear form on $\mathcal{NC}(A)$ by restriction of its domain.

\begin{defn}
A linear form $\Phi \in \overline{{\mathcal L_{NC}}(A)}$ is called a character if it is unital, $\Phi(\un)=1$, and multiplicative, i.e., for all $a,b \in T(\mathcal{NC}(A))$
$$
	\Phi(a|b)=\Phi(a)\Phi(b).
$$
A linear form $\kappa\in \overline{{\mathcal L_{NC}}(A)}$ is called an infinitesimal character, if $\kappa(\un)=0$, and if for all $a,b\in T(\mathcal{NC}(A))$ 
$$
	\kappa(a|b)=0.
$$
\end{defn}

We write ${G}_{NC}(A)$, respectively ${G}(A)$ for the set of characters in $\overline{{\mathcal L_{NC}}(A)}$, resp. $\overline{{\mathcal L}_A}$. We write $g_{NC}(A)$, respectively $g(A)$ for the corresponding sets of infinitesimal characters.

\begin{lem}\label{lemMap}
Since the map $Sp$ from $T(T(A))$ to $T({\mathcal NC}(A))$ is induced multiplicatively by a map from $T(A)$ to ${\mathcal NC}(A)$, the map $Sp^\ast$ restricts to maps from ${G}_{NC}(A)$ to ${G}(A)$ and $g_{NC}(A)$ to $g(A)$.
\end{lem}

\begin{thm}\label{thm:Gg}
There exists a natural bijection between $G_{NC}(A)$, the set of characters, and $g_{NC}(A)$, the set of infinitesimal characters on $\overline T(\mathcal{NC}(A))$. More precisely, for $\Phi \in G_{NC}(A), \exists ! \kappa  \in g_{NC}(A)$ such that 
$$
	\Phi=e+\kappa\prec \Phi = \exp^{\prec}( \kappa),
$$
and conversely, for $\kappa\in g_{NC}(A)$
$$
	\Phi:=\exp^{\prec}( \kappa)
$$
is a character.
This isomorphim commutes with $Sp^\ast$, in the sense that, given $\Phi$ and $\kappa$ as in the first part of the Theorem, $$Sp^\ast(\Phi)=e+Sp^\ast(\kappa)\prec Sp^\ast(\Phi) = \exp^{\prec}( Sp^\ast(\kappa)).$$
\end{thm}

The proof of the first part of the theorem is the same as the one of Theorem 7 for $G(A)$ and $g(A)$ in \cite{EPnc} and is omitted. The second part follows from Theorem \ref{thm:splitHomdual} and Lemma \ref{lemMap}.

\ \par

 Notice that elements in  ${G}_{NC}(A)$ and $g_{NC}(A)$ are entirely characterized by their restrictions to $\mathcal{NC}(A)$; similarly elements in ${G}(A)$ and $g(A)$ are characterized by their restrictions to $T(A)$. It follows that any section $\sigma$ of the map $Sp^\ast$ from $Lin(\mathcal{NC}(A),\mathbb{K})$ to $Lin(T(A),\mathbb{K})$ induces a surjection from  ${G}_{NC}(A)$ (resp. $g_{NC}(A)$) to ${G}(A)$ (resp. $g(A)$). 

In view of Theorem \ref{thm:Gg}, one can therefore use such a section to lift the equations $\Phi=e+\kappa\prec\Phi$ relating free moments and free cumulants in Theorems \ref{tim:freeprob} and \ref{gencumulant} to non-crossing partitions, that is, to $Lin(\mathcal{NC}(A),\mathbb{K})$. This process can be achieved in two different ways: either through the surjection from  ${G}_{NC}(A)$  to ${G}(A)$, or through the one from $g_{NC}(A)$ to $g(A)$. The second process happens to be the one best appropriated for our purposes.

We introduce for that reason a ``standard section'' that will lead to a new presentation of the classical M\"obius-inversion type relations between free moments and free cumulants. We however point out that \it other choices \rm of sections are possible that would lead to other lifts to non-crossing partitions of the relations between free moments and free cumulants.

\begin{defn}
Let $\kappa$ be a unital map from $T(A)$ to $\mathbb{K}$. We call standard section of $\kappa$ the linear form $sd(\kappa)$ on $\mathcal{NC}(A)$ defined by
$$
	sd(\kappa)(L \otimes a_1\cdots a_n):=\kappa(a_1\cdots a_n)
$$
if $L$ is the trivial non-crossing partition ($\hat{1}_n=L =[n]$), and zero else.
\end{defn}

For an arbitrary non-crossing partition $L=\{L_1,\ldots ,L_k\}$ of $[n]$, we write $\kappa^L(a_1\cdots a_n):=\prod\limits_{i=1}^k\kappa(a_{L_i})$.

\begin{prop}\label{keyrell}
The solution $\Psi$ of the equation $\Psi=e+sd(\kappa)\prec\Psi$ in $Lin(\mathcal{NC}(A),\mathbb{K})$ satisfies the identity:
\begin{equation}
\label{principEq}
	\Psi(Sp(a_1\cdots a_n))=\sum\limits_{L\in NC_n}\kappa^L(a_1\cdots a_n)=\phi(a_1\cdots a_n),
\end{equation}
where we recognize the standard relation between free moments and non-crossing cumulants (\ref{stdrel}).
More precisely, we have, for an arbitrary $L\in NC_n$,
\begin{equation}
\label{fundEqn}
\Psi(L\otimes a_1\cdots a_n)=\kappa^L(a_1\cdots a_n).
\end{equation}
\end{prop}

The proof goes over by induction on $[n]$. Let us assume that identity (\ref{fundEqn}) holds for non-crossing partitions of $[p]$, $p<n$. We get:
$$
	D:=\Psi(L\otimes a_1\cdots a_n)
	    = \sum\limits_{Q\!\coprod\limits_{\rm{adm}}\!P=L}st(\kappa)
	    \big(st(Q)\otimes a_Q\big)\Psi\big(st(J_1^{L,Q})\otimes
		 		a_{J_1^{L,Q}}|\cdots |st(J_{k(L,Q}^{L,Q})\otimes a_{J_{k(L,Q)}^{L,Q}}\big).
$$
However, since $sd(\kappa)$ vanishes on all non-crossing partitions excepted the trivial ones, terms on the right hand-side vanish except when $Q$ is the component of $L$ containing 1, and the expression reduces finally to
$$
	D=sd(\kappa)(st(L_1)\otimes a_{L_1})\Psi(st(L_2)\otimes a_{L_2}| \cdots | st(L_k)\otimes a_{L_k})
$$
$$
	D=\kappa( a_{L_1})\Psi(st(L_2)\otimes a_{L_2}| \cdots | st(L_k)\otimes a_{L_k})
$$
where $L=\{L_1,\ldots ,L_k\}$. From the induction hypothesis, we get the expected identity:
$$
	D=\prod\limits_{i=1}^k\kappa(a_{L_i}).
$$


\end{document}